\documentclass[conference]{IEEEtran}
\IEEEoverridecommandlockouts
\usepackage{cite}
\usepackage{amsmath,amssymb,amsfonts}
\usepackage{algorithmic}
\usepackage{graphicx}
\usepackage{textcomp}
\usepackage{xcolor}
\def\BibTeX{{\rm B\kern-.05em{\sc i\kern-.025em b}\kern-.08em
    T\kern-.1667em\lower.7ex\hbox{E}\kern-.125emX}}
\begin{document}

\title{Resilience-Oriented Parametric Insurance Design for Power Systems Under Extreme Weather
% \thanks{Funding information to be added if applicable.}
}

\author{
\IEEEauthorblockN{Jing Huang, Yawen Ma, Jiale Guo, and Chenjia Gu*}
\IEEEauthorblockA{
\textit{College of Electrical Engineering, Sichuan University}\\
Sichuan, China\\
huangjing123@scu.edu.cn, 2025141440416@stu.scu.edu.cn,\\
guojiale1030@foxmail.com, gcj0629@scu.edu.cn
}
\thanks{This work was supported by the China Postdoctoral Science Foundation under Grant 2025M770483.}
}

\maketitle

\begin{abstract}
Extreme weather leaves power systems exposed to residual outage risk even after physical resilience investments. Parametric insurance can provide pre-agreed contingent liquidity, but its physical value depends on how trigger thresholds and payout levels are designed. This paper proposes a resilience-oriented parametric insurance framework that couples a three-tier wind-index contract with post-event network restoration. Insurance payout expands the budget available to activate emergency resources, so the contract changes the physical restoration feasible set rather than merely offsetting accounting losses. Trigger thresholds and payout levels are jointly designed to balance actuarial premium, expected post-event system cost, and the conditional value-at-risk (CVaR) of scenario energy not supplied (ENS). A response-library method precomputes the restoration mixed-integer linear program for each scenario-payout pair and then evaluates admissible contracts efficiently. On the IEEE RTS-24 with 80 extreme-wind scenarios, the optimized contract reduces expected EENS and $\mathrm{CVaR}_{0.90}$ of ENS by 21.1\% and 21.4\%, respectively, relative to no insurance, while requiring 48.8\% less premium than a fixed parametric contract with comparable resilience. The results show that insurance design should target the nonlinear liquidity-to-resilience response rather than loss compensation alone.
\end{abstract}

\begin{IEEEkeywords}
parametric insurance, power system resilience, extreme weather, emergency restoration, CVaR, risk management.
\end{IEEEkeywords}

\section{Introduction}
Extreme weather can simultaneously damage multiple network components, isolate load pockets, and create prolonged service interruptions. Such low-probability, high-impact events have motivated extensive work on fragility modeling, probabilistic resilience assessment, hardening, and operational adaptation \cite{Panteli2017}. More broadly, climate-adaptive operation and planning must account for both acute extreme-event disruptions and longer-term shifts in renewable resources and planning boundaries under nonstationary climate conditions \cite{Gu2026Climate}. Nevertheless, physical measures cannot economically eliminate all residual risk \cite{Gu2026Planning}. Financial instruments therefore provide a complementary mechanism for transferring or pooling losses that remain after preventive investments.

Insurance has recently been incorporated into power-system risk management from several perspectives. Catastrophe insurance has been used to allocate distribution-system disaster risk \cite{Sun2023}, and resilience planning models have linked insurance decisions with infrastructure investment \cite{Hu2024}. Insurance-inspired mechanisms have also been proposed to coordinate distributed resilience investment \cite{Billimoria2023,Huang2026Insurance}, while actuarial frameworks are emerging for transmission-system extreme-weather risk \cite{Zhao2026}. These studies demonstrate the relevance of insurance to power-system resilience, but insurance is often represented mainly as financial compensation, premium allocation, or an investment signal. The direct pathway from an insurance payout to the physically feasible post-event restoration actions is less explicit.

Parametric insurance is attractive for this purpose because payment is determined by an observable index rather than detailed ex-post loss adjustment. A pre-agreed trigger can therefore support rapid access to contingent liquidity. The drawback is \emph{basis risk}: an index-based payout may not match the realized network loss because damage depends on spatial exposure, component fragility, and network topology \cite{Niakh2025}. Hence, simply increasing coverage does not necessarily maximize resilience. Once emergency resource capacity or deliverability becomes binding, additional payout can raise the premium without further reducing unserved energy.

This paper develops a compact framework for designing parametric insurance specifically around this physical liquidity-to-resilience relationship. The main contributions are twofold. First, the scenario-dependent payout is embedded in a post-event emergency budget constraint, allowing insurance to alter the feasible activation and dispatch of emergency resources. Second, wind-speed trigger thresholds and tiered payout levels are jointly designed to minimize a risk-adjusted objective comprising premium, expected restoration cost, and tail ENS. A response-library solution exploits the fact that restoration depends on a contract only through its realized payout, which permits efficient evaluation of a finite set of candidate contracts.

\begin{figure*}[t]
\centering
\includegraphics[width=0.96\textwidth]{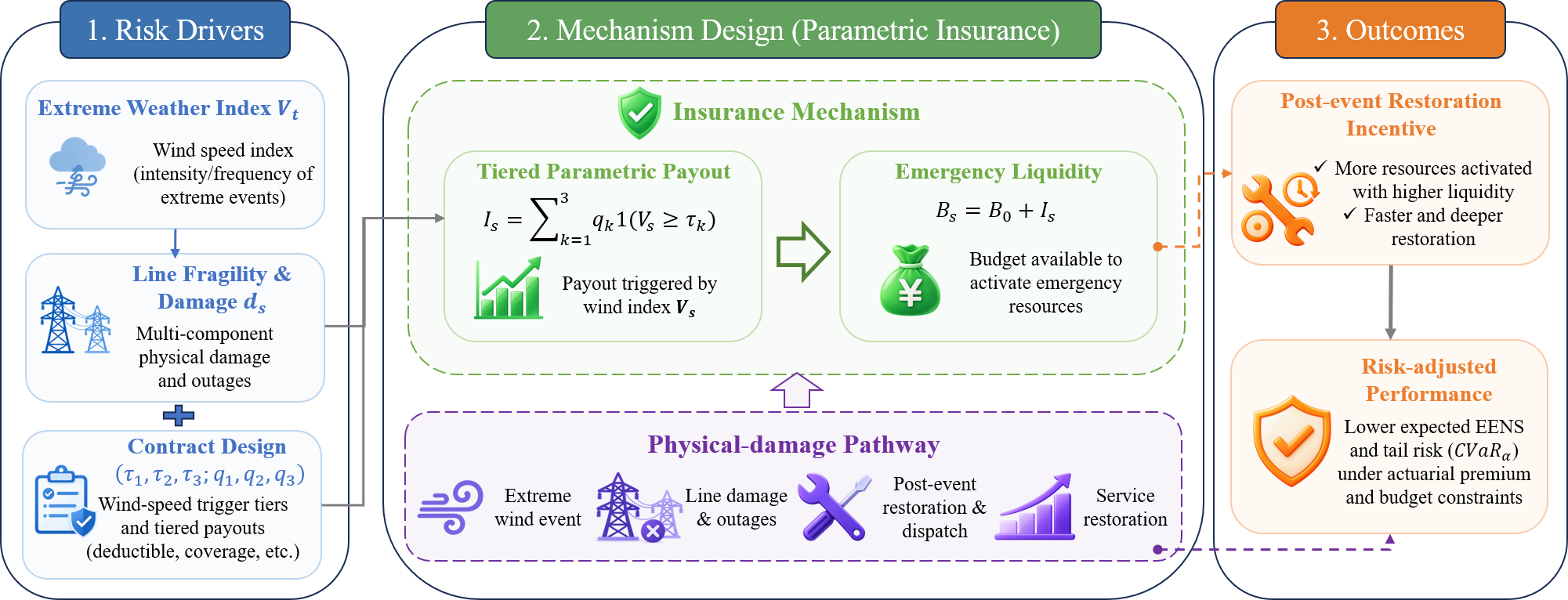}
\caption{Resilience-oriented parametric insurance mechanism. }
\label{fig:framework}
\end{figure*}

\section{Parametric Insurance-Enabled Restoration Model}
\subsection{Extreme-Weather Scenarios and Parametric Contract}
Let $s\in\mathcal S$ denote an extreme-weather scenario with probability $p_s$, and let $V_s$ be the regional peak-wind index. Spatially heterogeneous line wind is represented by
\begin{equation}
    v_{l,s}=V_s\alpha_l\xi_{l,s},
\end{equation}
where $\alpha_l$ is a line exposure factor and $\xi_{l,s}$ is a local perturbation. A logistic fragility model gives the conditional failure probability
\begin{equation}
    \rho_{l,s}=\frac{1}{1+\exp[-\beta(v_{l,s}-v_l^{50})]},
\end{equation}
where $v_l^{50}$ is the median failure wind speed and $\beta$ controls the fragility slope. Sampling $\rho_{l,s}$ yields the line-availability parameter $a_{l,s}\in\{0,1\}$.

A three-tier parametric contract is characterized by increasing trigger thresholds $\boldsymbol\tau=[\tau_1,\tau_2,\tau_3]$ and payouts $\mathbf q=[q_1,q_2,q_3]$:
\begin{equation}
I_s(\boldsymbol\tau,\mathbf q)=
\begin{cases}
0, & V_s<\tau_1,\\
q_1, & \tau_1\le V_s<\tau_2,\\
q_2, & \tau_2\le V_s<\tau_3,\\
q_3, & V_s\ge\tau_3.
\end{cases}
\label{eq:payout}
\end{equation}
The actuarial premium is modeled by an expected-value principle,
\begin{equation}
    \Pi=(1+\lambda)\sum_{s\in\mathcal S}p_s I_s,
\label{eq:premium}
\end{equation}
where $\lambda$ is the loading factor. The claim is assumed to be determined at the beginning of the restoration horizon; operationally, $I_s$ can be interpreted as payout or immediately available bridge liquidity backed by the triggered contract.

\subsection{Post-Event Network Restoration}
For a given scenario and payout, the system operator minimizes post-event cost over time periods $t\in\mathcal T$:
\begin{align}
C_s^*(I_s)=\min\ &\sum_{g,t}c_g P_{g,t,s}
+\sum_e F_e y_{e,s}
+\sum_{e,t}c_e^{E}P_{e,t,s}^{E} 
\nonumber\\
&+\sum_{i,t}v_i P_{i,t,s}^{LS}. \label{eq:recourseobj}
\end{align}
Here $P_{g,t,s}$ is conventional generation, $P_{e,t,s}^{E}$ is emergency-resource output, $y_{e,s}$ is its binary activation variable, $P_{i,t,s}^{LS}$ is load shedding, and $v_i$ is the value of lost load (VOLL). Critical loads are assigned a larger $v_i$.

The key insurance-resilience coupling is
\begin{equation}
\sum_e F_e y_{e,s}+\sum_{e,t}c_e^{E}P_{e,t,s}^{E}
\le B_0+I_s,
\label{eq:liquidity}
\end{equation}
where $B_0$ is baseline emergency liquidity. Therefore, the payout does not simply reduce an accounting loss; it enlarges the feasible set of emergency actions. Emergency resource limits are
\begin{equation}
0\le P_{e,t,s}^{E}\le \bar P_e^{E}y_{e,s}.
\end{equation}

A DC network representation is used. Line availability is imposed directly through
\begin{align}
-\bar f_l a_{l,s}\le f_{l,t,s}\le \bar f_l a_{l,s}, \label{eq:linecap}\\
-M_l(1-a_{l,s})\le f_{l,t,s}-b_l(\theta_{i,t,s}-\theta_{j,t,s})
\le M_l(1-a_{l,s}), \label{eq:dcflow}
\end{align}
for line $l=(i,j)$, where $M_l$ relaxes the angle--flow relation when the line fails. Thus $a_{l,s}=0$ forces $f_{l,t,s}=0$, while an intact line follows the DC power-flow relation. Nodal balance is
\begin{equation}
\sum_{g\in\mathcal G_i}P_{g,t,s}+\sum_{e\in\mathcal E_i}P_{e,t,s}^{E}
-P_{i,t}^{D}+P_{i,t,s}^{LS}=\sum_{l}A_{il}f_{l,t,s}. \label{eq:balance}
\end{equation}
The remaining operating bounds are
\begin{equation}
0\le P_{g,t,s}\le\bar P_g,\qquad
0\le P_{i,t,s}^{LS}\le P_{i,t}^{D}. \label{eq:opbounds}
\end{equation}
The zero lower bound on conventional generation is intentional: a generator located in a post-event island can be shut down rather than creating an artificial infeasibility through its pre-contingency minimum-output limit. The scenario energy not supplied (ENS) is
\begin{equation}
    E_s=\sum_{i,t}P_{i,t,s}^{LS}\Delta t. \label{eq:ens}
\end{equation}

\subsection{Risk-Aware Contract Design}
The contract is selected from ordered candidate sets $\mathcal T_\tau$ and $\mathcal Q$, subject to $\tau_1<\tau_2<\tau_3$ and $0<q_1<q_2<q_3$. We minimize
\begin{equation}
\min_{\boldsymbol\tau,\mathbf q}\ \Pi+\sum_s p_s C_s^*(I_s)
+\omega \bar v\,\mathrm{CVaR}_{\alpha}(E_s),
\label{eq:design}
\end{equation}
where $\omega$ is a tail-risk weight and $\bar v$ monetizes ENS in the risk term. The expected value $\sum_s p_s E_s$ is EENS. The discrete-scenario CVaR term is represented using auxiliary variables $\eta$ and $z_s$ \cite{Rockafellar2000}:
\begin{align}
\mathrm{CVaR}_{\alpha}(E)&=\eta+\frac{1}{1-\alpha}\sum_s p_s z_s, \label{eq:cvar}\\
z_s&\ge E_s-\eta,\qquad z_s\ge0. \label{eq:cvaraux}
\end{align}
The expected-cost term already internalizes ordinary load shedding through VOLL, while the CVaR term provides an explicit preference against severe tail-ENS scenarios.

For diagnostic purposes, an indemnity-reference payout is defined from the no-insurance physical loss,
\begin{equation}
I_s^{\rm ref}=\min\{\kappa\bar v E_s^0,\bar I\}, \label{eq:indref}
\end{equation}
where $\kappa$ is a reference coverage fraction and $\bar I$ is the payout cap. Basis risk is then measured by
\begin{equation}
BR=\sum_s p_s\left|I_s-I_s^{\rm ref}\right|. \label{eq:basis}
\end{equation}
This metric is not included in \eqref{eq:design}: the contract is deliberately designed for physically useful restoration liquidity rather than dollar-for-dollar loss replication.

Directly nesting a mixed-integer restoration model inside every contract evaluation is unnecessary. For each scenario $s$ and payout level $q\in\mathcal Q$, the restoration problem is solved once and the response tuple $\mathcal R_{s,q}=\{C_{s,q}^*,E_{s,q},E_{s,q}^{\rm crit}\}$ is stored. A candidate contract maps $V_s$ to one payout level and therefore to one stored response. With $|\mathcal S|=80$ and seven payout levels, the response library requires only 560 restoration MILPs; all admissible contracts are then evaluated by lookup. This decomposition also exposes the scenario-specific marginal physical value of liquidity, $E_{s,q_a}-E_{s,q_b}$, which is used in the case-study interpretation.

\section{Case Study}
\subsection{Test System and Settings}
The framework is tested on the IEEE RTS-24 \cite{RTS1979}, using the MATPOWER data format \cite{Zimmerman2011}. The study uses 80 equally likely, reproducible extreme-wind scenarios and a six-hour post-event horizon. Six dispatchable emergency resources, each rated at 120~MW, are placed at buses 3, 10, 13, 15, 18, and 20. The five largest-load buses are treated as critical loads. Line outages remain fixed over the short restoration horizon so that the study isolates the interaction between contingent liquidity and emergency dispatch rather than crew routing or repair sequencing.

\begin{table}[t]
\caption{Main Case-Study Parameters}
\label{tab:parameters}
\centering
\footnotesize
\setlength{\tabcolsep}{3.5pt}
\begin{tabular}{p{0.52\columnwidth}p{0.39\columnwidth}}
\hline
Parameter & Value \\
\hline
Scenarios / horizon & 80 / 6 h \\
Peak wind $V_s$ & truncated $\mathcal N(44,7.5^2)$, $[30,62]$ m/s \\
Line exposure $\alpha_l$ & $U[0.85,1.05]$ \\
Fragility $\beta$, nominal $v^{50}$ & 0.24 / 54 m/s \\
Baseline liquidity $B_0$ & 250 kUSD \\
Emergency resources & $6\times120$ MW \\
Activation / variable cost & 40 kUSD / 0.35 kUSD/MWh \\
Normal / critical VOLL & 10 / 30 kUSD/MWh \\
Threshold candidates & 36:4:56 m/s \\
Payout grid & 0:200:1200 kUSD \\
Loading $\lambda$ / CVaR $\alpha$ & 0.15 / 0.90 \\
Risk weight $\omega$ & 0.25 \\
\hline
\end{tabular}
\end{table}

The event generator yields an average of 6.9 failed lines per scenario, with a 90th-percentile outage count of approximately 16. The finite candidate sets in Table~\ref{tab:parameters} produce $\binom{6}{3}\binom{6}{3}=400$ admissible contracts. A fixed contract with thresholds $[40,48,52]$~m/s and payouts $[200,600,1200]$~kUSD is used as a hand-designed benchmark. All restoration problems are linearized mixed-integer models; the supplied implementation uses YALMIP/Gurobi.

\subsection{Optimized Contract and Resilience Improvement}
The optimized contract has trigger thresholds
\begin{equation*}
\boldsymbol\tau^*=[48,52,56]\ \mathrm{m/s}
\end{equation*}
and payout levels
\begin{equation*}
\mathbf q^*=[400,600,1000]\ \mathrm{kUSD}.
\end{equation*}
Its trigger rate is 35\%, giving an expected payout of 212.5~kUSD and a loaded premium of 244.4~kUSD. Table~\ref{tab:results} compares the three cases.

\begin{table}[t]
\caption{Comparison of Insurance Designs}
\label{tab:results}
\centering
\scriptsize
\setlength{\tabcolsep}{2.2pt}
\begin{tabular}{lrrrr}
\hline
Case & Premium & Exp. EENS & CVaR ENS & Objective \\
 & (kUSD) & (MWh) & (MWh) & (kUSD)\\
\hline
No insurance & 0.0 & 765.8 & 5044.6 & 23487.2\\
Fixed param. & 477.3 & 602.7 & 3966.1 & 19212.8\\
Optimized param. & 244.4 & 604.2 & 3966.1 & 18994.7\\
\hline
\end{tabular}
\end{table}

\begin{figure}[t]
\centering
\includegraphics[width=0.98\columnwidth]{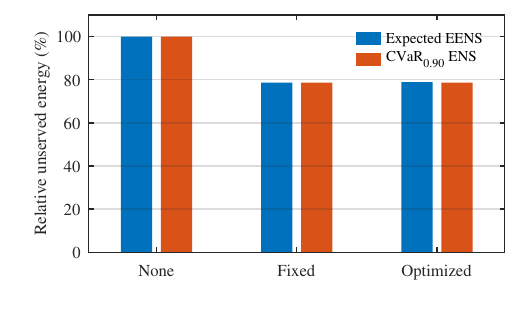}
\caption{Expected and tail unserved energy normalized by the no-insurance case. Exact values are reported in Table~\ref{tab:results}.}
\label{fig:resilience}
\end{figure}

Relative to no insurance, the optimized contract reduces expected EENS from 765.8 to 604.2~MWh (21.1\%) and $\mathrm{CVaR}_{0.90}$ of ENS from 5044.6 to 3966.1~MWh (21.4\%). Critical-load EENS falls from 146.9 to 122.4~MWh (16.7\%). Fig.~\ref{fig:resilience} shows that the fixed and optimized contracts obtain nearly identical average and tail resilience. Their financial requirements, however, differ markedly.

The optimized trigger structure is selective rather than uniformly generous. Among the 80 scenarios, 52 receive no payout, while 11, 11, and 6 scenarios receive 400, 600, and 1000~kUSD, respectively. Of the 28 triggered scenarios, 15 exhibit an actual EENS reduction. Average EENS reductions within the three payout tiers are 225.7, 316.2, and 1160.8~MWh, respectively. The remaining triggered scenarios receive liquidity but show no physical improvement because emergency-resource deliverability, rather than the budget, is binding. This distinction explains why the optimized contract can use a substantially lower premium without sacrificing tail resilience.

\subsection{Scenario-Level Liquidity Response}
The aggregate metrics conceal an important source of heterogeneity: the same amount of contingent liquidity can have very different physical value across outage topologies. Table~\ref{tab:scenarioresponse} reports two representative scenarios from the response library. In Scenario A, 16 failed lines create substantial but partially recoverable load isolation. Increasing available insurance liquidity from 0 to 400~kUSD allows two additional emergency resources to be activated and raises emergency energy from 501.1 to 1387.7~MWh. EENS consequently falls from 2314.8 to 1428.2~MWh. Raising the payout further to 1200~kUSD produces no additional reduction, indicating that the useful emergency actions are already exhausted at 400~kUSD.

\begin{table}[t]
\caption{Representative Scenario Response to Additional Liquidity}
\label{tab:scenarioresponse}
\centering
\scriptsize
\setlength{\tabcolsep}{2.0pt}
\begin{tabular}{lrrrrr}
\hline
Scenario & $V_s$ & Failed & \multicolumn{3}{c}{EENS (MWh)}\\
 & (m/s) & lines & $I=0$ & $I=400$ & $I=1200$\\
\hline
A & 51.63 & 16 & 2314.8 & 1428.2 & 1428.2\\
B & 49.62 & 11 & 1096.1 & 1096.1 & 1096.1\\
\hline
\end{tabular}
\end{table}

Scenario B illustrates the complementary case. Although 11 lines fail and EENS reaches 1096.1~MWh, additional payout does not change the solution. Only about 10~MWh of emergency energy is deliverable even when the liquidity limit is relaxed, because the relevant load pocket cannot be effectively reached from the available emergency-resource locations. Thus, insurance is most valuable when financial scarcity is the active recovery bottleneck; it cannot substitute for network connectivity, resource siting, or repair capability. In the response library, EENS is nonincreasing in payout for every scenario, but it is often piecewise flat because additional dollars do not continuously create additional deliverable emergency power. This scenario dependence is precisely why contract design based only on expected monetary loss can over-insure some events while under-serving others.

\subsection{Premium-Resilience Trade-off and Basis Risk}
\begin{figure}[t]
\centering
\includegraphics[width=0.98\columnwidth]{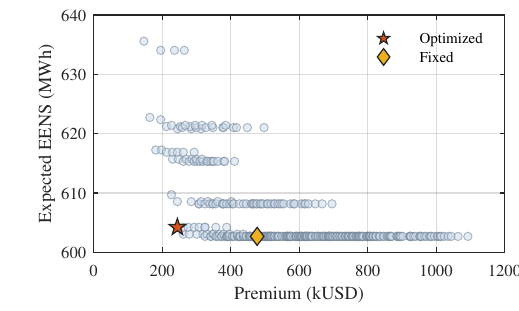}
\caption{Premium-resilience trade-off across 400 admissible contracts.}
\label{fig:tradeoff}
\end{figure}

Fig.~\ref{fig:tradeoff} plots all 400 candidate contracts in the premium-EENS plane. The discrete bands reflect the finite payout levels and the nonlinear restoration response. The hand-designed fixed contract attains slightly lower expected EENS (602.7 versus 604.2~MWh), but its premium is 477.3~kUSD. The optimized design reaches essentially the same physical resilience with a 244.4~kUSD premium, a 48.8\% reduction. This occurs because additional liquidity becomes ineffective once the useful emergency resources at electrically relevant locations are already activated or network transfer constraints become binding. The optimizer therefore avoids paying for coverage that provides little marginal physical benefit.

\begin{figure}[t]
\centering
\includegraphics[width=0.98\columnwidth]{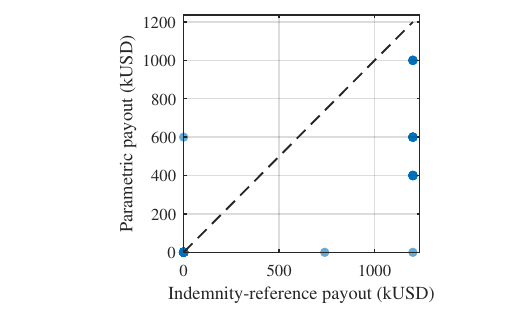}
\caption{Basis risk of the optimized contract relative to the indemnity-reference payout. The dashed line denotes payout equality.}
\label{fig:basis}
\end{figure}

To diagnose basis risk, an idealized indemnity reference is defined as 60\% of the no-insurance VOLL-based unserved-energy loss, capped at the maximum payout of 1200~kUSD. This reference is used only for evaluation and not in the design objective. Fig.~\ref{fig:basis} compares it with the optimized parametric payout. The mean absolute payout mismatch is 231.7~kUSD. The optimized contract mainly underpays the indemnity reference because it is designed to purchase the liquidity that is operationally useful for restoration, rather than to replicate financial loss dollar-for-dollar. This distinction is central: a resilience-oriented parametric contract need not be an accurate loss estimator if its purpose is to fund effective recovery actions. Nevertheless, large basis risk can leave some high-loss events underfunded, motivating richer spatial or multi-index triggers in future work.

\subsection{Premium-Constrained Design and Marginal Value of Coverage}
\begin{figure}[t]
\centering
\includegraphics[width=0.98\columnwidth]{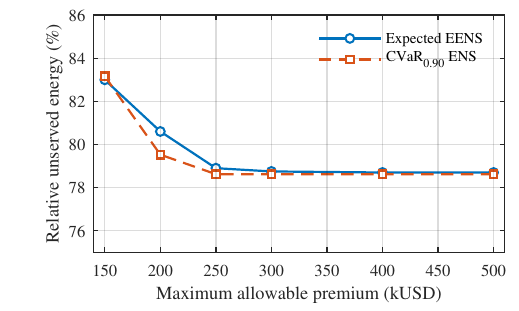}
\caption{Best attainable expected and tail unserved energy under alternative premium caps, normalized by the no-insurance case.}
\label{fig:premiumcap}
\end{figure}

The full contract set also permits a budget-oriented interpretation that is useful for utility decision makers. For a premium cap $\bar\Pi$, define the best attainable physical performance as
\begin{equation}
E^{\min}(\bar\Pi)=\min_{(\boldsymbol\tau,\mathbf q):\,\Pi\le\bar\Pi}
\sum_s p_s E_s.
\label{eq:premiumcap}
\end{equation}
Fig.~\ref{fig:premiumcap} evaluates this envelope using the same 400 contracts, so no additional restoration optimization is required. With a 150~kUSD cap, the selected contract costs 146.6~kUSD and lowers expected EENS to 635.6~MWh, a 17.0\% reduction from the no-insurance case. Raising the cap to 250~kUSD allows the optimized contract in Table~\ref{tab:results}, reducing expected EENS and tail ENS by 21.1\% and 21.4\%, respectively. Beyond about 300~kUSD, however, the improvement is negligible: increasing the feasible premium to 400~kUSD changes the best EENS only from 603.1 to 602.7~MWh, while tail ENS is unchanged. This confirms that the relevant quantity for contract design is not payout size itself but the marginal physical value of additional liquidity.

The saturation region has two operational interpretations. First, once the emergency resources that can effectively serve disconnected load pockets are activated, a larger budget cannot create additional deliverable power. Second, surviving transmission corridors can become the binding limitation, so extra emergency generation may be stranded behind network constraints. Consequently, the optimal insurance design should be coordinated with the location and capability of emergency resources rather than priced solely from aggregate expected loss. This result also suggests a natural planning extension in which mobile-resource siting, hardening, and insurance parameters are co-designed under evolving climate risk.

\section{Conclusion}
This paper proposed a resilience-oriented parametric insurance design in which a weather-index payout directly enlarges the emergency budget available for post-event power-system restoration. Joint design of wind-speed triggers and tiered payouts balances premium cost, expected restoration performance, and tail ENS. The IEEE RTS-24 study shows that the optimized contract can reduce both expected and tail unserved energy by about 21\% while using roughly half the premium of a more generous fixed contract with comparable resilience. The results also reveal a saturation effect: beyond the liquidity needed to activate operationally useful emergency resources, additional coverage may have little physical value. Future work will replace stylized fragility scenarios with empirically calibrated hazard models and incorporate repair crews, settlement delays, insurer solvency, and spatially richer triggers.


\begin{thebibliography}{12}
\bibitem{Panteli2017}
M. Panteli, C. Pickering, S. Wilkinson, R. Dawson, and P. Mancarella, ``Power system resilience to extreme weather: Fragility modeling, probabilistic impact assessment, and adaptation measures,'' \emph{IEEE Trans. Power Syst.}, vol. 32, no. 5, pp. 3747--3757, Sep. 2017.


\bibitem{Gu2026Climate}
C. Gu, J. Ruan, X. Yang, J. Huang, Y. Qiu, J. Wang, S. Chen, G. Liang, Z. Xu, H. Su, T. T. Lie, and J. Song, ``Advancing climate-adaptive operation and planning for renewable-rich energy systems,'' \emph{Renew. Sustain. Energy Rev.}, vol. 242, Art. no. 117336, 2026, doi: 10.1016/j.rser.2026.117336.

\bibitem{Gu2026Planning}
C. Gu, J. Ruan, Y. Qiu, et al., ``Toward climate-adaptive low-carbon power system planning: A multistage stochastic framework considering climate uncertainties,'' \emph{IEEE Trans. Ind. Informat.}, vol. 22, no. 5, pp. 3938--3949, 2026.

\bibitem{Sun2023}
S. Sun, G. Li, Y. Bian, Z. Bie, and Q. Hu, ``Catastrophe risk management for electric power distribution systems: An insurance approach,'' \emph{CSEE J. Power Energy Syst.}, vol. 9, no. 1, pp. 393--410, Jan. 2023.

\bibitem{Hu2024}
Q. Hu, G. Li, S. Sun, and Z. Bie, ``Incorporating catastrophe insurance in power distribution systems investment and planning for resilience enhancement,'' \emph{Int. J. Electr. Power Energy Syst.}, vol. 155, Art. no. 109438, Jan. 2024.

\bibitem{Billimoria2023}
F. Billimoria, F. Fele, I. Savelli, T. Morstyn, and M. McCulloch, ``An insurance paradigm for improving power system resilience via distributed investment,'' \emph{IEEE Trans. Energy Markets, Policy Regul.}, vol. 1, no. 4, pp. 499--511, Dec. 2023.

\bibitem{Huang2026Insurance}
J. Huang, K. Liao, C. Gu, et al., ``Enhancing power distribution system resilience through insurance mechanisms: An insurer optimization approach with V2G aggregator participation,'' \emph{IEEE Trans. Smart Grid}, 2026.

\bibitem{Zhao2026}
Y. Zhao, Y. Ding, X. Sun, W. Shi, Z. Hu, Z. Xu, C. Chen, and Z. Bie, ``A multi-stakeholder actuarial risk management framework for transmission systems against extreme weather events,'' \emph{IEEE Trans. Energy Markets, Policy Regul.}, vol. 4, no. 1, pp. 78--92, 2026.

\bibitem{Niakh2025}
F. Niakh, A. Bassi\`ere, M. Denuit, and C. Y. Robert, ``Peer-to-peer basis risk management for renewable production parametric insurance,'' \emph{Ann. Oper. Res.}, 2025, doi: 10.1007/s10479-025-06987-w.

\bibitem{Rockafellar2000}
R. T. Rockafellar and S. Uryasev, ``Optimization of conditional value-at-risk,'' \emph{J. Risk}, vol. 2, no. 3, pp. 21--41, 2000.

\bibitem{RTS1979}
IEEE Reliability Test System Task Force, ``IEEE reliability test system,'' \emph{IEEE Trans. Power App. Syst.}, vol. PAS-98, no. 6, pp. 2047--2054, Nov./Dec. 1979.

\bibitem{Zimmerman2011}
R. D. Zimmerman, C. E. Murillo-S\'anchez, and R. J. Thomas, ``MATPOWER: Steady-state operations, planning, and analysis tools for power systems research and education,'' \emph{IEEE Trans. Power Syst.}, vol. 26, no. 1, pp. 12--19, Feb. 2011.
\end{thebibliography}
\end{document}